\documentclass[%
times
,preprint
,1p
,authoryear
]{elsarticle}

\usepackage[utf8]{inputenc}
\usepackage{amsfonts,amssymb,amsmath,amsthm,stmaryrd,enumerate}
\usepackage{lmodern,bbold}  
\usepackage{mathtools,etoolbox,xcolor}
\usepackage[mathcal]{euscript}
\usepackage{verbatim}
\usepackage{hhline}
\usepackage{caption}

\graphicspath{ {./images/} }

\usepackage[citecolor=teal,linkcolor=teal,
    urlcolor=teal,
    colorlinks]{hyperref}

\newtheorem{theorem}{Theorem}
\newtheorem{lemma}{Lemma}
\newtheorem{corollary}{Corollary}
\newtheorem{proposition}{Proposition}
\newdefinition{remark}{Remark}
\newtheorem*{conj}{Conjecture}

\makeatletter
\def\ps@pprintTitle{%
  \let\@oddhead\@empty
  \let\@evenhead\@empty
  \let\@oddfoot\@empty
  \let\@evenfoot\@oddfoot
}
\makeatother

\newcommand{\Var}{\operatorname{Var}}

\newcommand{\eqinlaw}{\stackrel{d}{=}}

\DeclarePairedDelimiter{\abs}{\lvert}{\rvert}
\DeclarePairedDelimiter{\zjel}{(}{)}
\DeclarePairedDelimiterX\zfrac[2]{(}{)}{%
  \frac{#1}{#2}%
}

\DeclarePairedDelimiter\br\lbrack\rbrack
\DeclarePairedDelimiter\event\lbrace\rbrace
\DeclarePairedDelimiterX\set[2]\{\}{%
  #1\nonscript\::\allowbreak\nonscript\:\mathopen{}#2%
}

\def\toby#1#2{\,{\buildrel #2\over#1}\,}

\newcommand\pto[1][\to]{\toby{#1}{p}}

\let\cbr\event

\let\phi\varphi
\let\theta\vartheta

\makeatletter
  \def\isparam{\@ifnextchar{\bgroup}}%
\makeatother

\let\texexp\exp
\def\exp{\texexp\isparam{\cbr*}{}}

\let\PEfont\mathbb
\providecommand\given{}

\def\redefgiven#1{%
  \renewcommand\given{%
    \nonscript\:%
    #1\vert%
    \allowbreak%
    \nonscript\:%
    \mathopen{}%
  }%
}
\DeclarePairedDelimiterX\PEzjel[1](){\redefgiven{\delimsize}#1}
\newcommand{\PE}[1][]{%
  \PEfont\Pe%
  \ifblank{#1}{}{_{#1}}%
  \isparam{\PEzjel*}{}%
}
\newcommand{\E}{\def\Pe{E}\PE}
\renewcommand{\P}{\def\Pe{P}\PE}

\makeatletter

\def\ifbrace{\@ifnextchar{\bgroup}}

\DeclareListParser{\dolist}{}

\def\stripdef@end#1{\egroup\def#1}

\def\stripdef#1#2{%
  \bgroup
  \escapechar=-1\relax
  \expandafter\stripdef@end\csname#1\string#2\endcsname
}
\makeatother

\def\defprefix#1#2{
    \def\do##1{\stripdef{#1}{##1}{{#2{##1}}}}%
    \ifbrace\dolist\do
}

\defprefix{t}{\tilde}{ABcXQYZmMS\mu\F\beta\tbeta\tau\Omega}
\defprefix{c}{\mathcal}{FGAHSCDLPBV}
\defprefix{}{\mathcal}{FGB}

\def\cadlag/{c\`adl\`ag}
\def\ito/{It\^o}
\def\sign{\operatorname{sign}}
\def\cov{\operatorname{cov}}
\def\real{\mathbb{R}}

\def\Isymb{{\mathbb 1}}
\def\substr#1{_{\zjel{#1}}}
\newcommand{\I}[1][]{%
  \Isymb\ifblank{#1}{\substr}{_{#1}}%
}

\def\dist{\operatorname{dist}}

\begin{document}
\begin{frontmatter}
  \title{Example of a Dirichlet process whose zero energy part has finite \texorpdfstring{$p$}{p}-variation} 
  \author[1]{Vilmos Prokaj}
  \ead{vilmos.prokaj@ttk.elte.hu}
  \author[2,1]{László Bondici}
  \ead{bondici@renyi.hu}
  \address[1]{ELTE Eötvös Loránd University, Budapest, Hungary.}
  \address[2]{Alfréd Rényi Institute of Mathematics, Budapest, Hungary.}
  \begin{abstract}
    Let $B^H$ be a fractional Brownian motion on $\real$ with Hurst parameter $H\in(0,1)$, 
    $F$ be its pathwise antiderivative with $F(0)=0$, %
    and let $B$ be a standard Brownian motion, independent of $B^H$. %
    We show that the zero energy part $A_t=F(B_t)-\int_0^t F'(B_s)dB_s$ of $F(B)$ has positive and 
    finite $p$-variation in a special sense for $p_0=\frac{2}{1+H}$. 
    We also present some simulation results about the zero energy part of a certain median process
    which suggest that its $4/3$-variation is positive and finite.
  \end{abstract}
\end{frontmatter}

\section{Introduction}

In a recent paper  (\cite{ProkajBondici2022}) we showed that a certain median process lacks the semimartingale property. 
This median process has a decomposition into a sum of a martingale and a process with zero quadratic variation. 
Such a process is called a Dirichlet process in \cite{follmer_1981}, and a process with zero quadratic variation
is said to have zero energy.
Since the median process above is not a semimartingale, the zero energy part of the decompostion can not be of finite total variation.
The proof provided 
in the above paper is indirect, that is it does not compute the total variation of the zero energy part. 
This paper tries to make the first steps for this computation. We consider a simpler case where the computation can be carried out
and we also present some simulation result for the above median process.
These simulation results suggest that if we compute the $4/3$-variation along a specially selected  sequence 
of partitions, then it has a finite limit. The exponent $4/3$ is the same as in \cite{MR1118450}.
They consider a process $X$ obtained from a Brownian motion $B$ with the formula $X_t=\int_0^t \I{B_s\leq B_t}ds$,  
that is $X$ is the amount of time spent so far  below the current value of $B$. 

Without going into too much details, short term increments 
of these processes are obtained roughly by substituting a Brownian motion $B$
into a continuously differentiable random function $F$ whose %
first derivative $f$ has non-zero and finite quadratic variation. The increment of the zero energy part of $F(B)$ 
during an upcrossing of the 
interval $(0,\delta)$ is $F(\delta)-F(0)-\int_0^\tau f(B_s)dB_s$ where $\tau$ is the exit time of $(-\delta,\delta)$ and 
$B$ is conditioned to $B_\tau=\delta$. By the mean value theorem there is a $\xi\in(0,\delta)$ such that 
$F(\delta)-F(0)=f(\xi)\delta$ and we get that 
the increment is roughly $\delta(f(\xi)-f(0))$. Here the stochastic integral is simply approximated with $f(0)(B_\tau-B_0)=f(0)\delta$. 

Similarly, the increment of the zero energy part during a downcrossing of $(0,\delta)$ is roughly 
$\delta(f(\delta)-f(\xi))$. Since the number of down and upcrossing upto time 
$t$ differs only by one and is proportional to $1/\delta$, 
we get that the contribution of crossings of $(0,\delta)$ is 
roughly proportional to $\frac1\delta\abs{\delta(f(\delta)-f(0))}^p$. Similar computation can be 
done for intervals of form $(r\delta,(r+1)\delta)$. 
The fact that the quadratic variation of $f$ is finite roughly means that 
for most of the increments $f((r+1)\delta)-f(r\delta)$ 
is of order $\delta^{1/2}$ so finally we arrive at the conclusion that in order to have finite 
$p$ variation it is needed that $\frac1\delta \delta^{3/2p}=\delta$, that is $p=4/3$. 

This argument is far from being rigorous. The aim of this paper to make this argument precise 
for the simplest case when $f$ is a fractional Brownian motion. The increments of a fractional Brownian motion 
with Hurst index $H$ over an interval of length $\delta$ 
is of size $\delta^{H}$, so the last step of the above heuristic computation is 
$\frac{1}{\delta}\delta^{(1+H)p}=\delta$, that is $p=2/(1+H)$.

\begin{theorem}\label{thm:1}
  Let $B^H$ be a fractional Brownian motion on $\real$ with Hurst
  parameter $H$ and let $B$ be a standard Brownian motion independent 
  of $B^H$. Denote $F$ an antiderivative of $B^H$ with $F(0)=0$ and %
  \begin{displaymath}
    A_t = F(B_t)-\int_0^t F'(B_s)dB_s.
  \end{displaymath}
  For $p_0=\frac{2}{1+H}$ %
  the $p_0$ variation of $A$ on any $[0,t]$ exists and equals to $ct$ with $c=\E{\abs{A_1}^{p_0}}$.

  Especially, for $p<p_0$ the $p$ variation is infinite almost surely 
  on any non degenerate interval, 
  while for $p>p_0$ it is identically zero.
\end{theorem}
In the previous claim  the $p$ variation $(V^{(p)}_t)_{t\geq0}$  is definied similarly to the quadratic variation,
that is, for any $t$ and any (deterministic) sequence of subdivisions $({t^{(n)}_i})$ of $[0,t]$ whose 
mesh goes to zero we have that
\begin{displaymath}
  \sum_i \abs*{A_{t^{(n)}_{i+1}}-A_{t^{(n)}_{i}}}^{p}\to V^{(p)}_t
  \quad
  \text{in probability}.
\end{displaymath}

We start with a somewhat simpler claim.
\begin{theorem}\label{thm:2}Using the notation of Theorem \ref{thm:1},
  for $\delta>0$ let 
  \begin{displaymath}
    \tau^{\delta}_0=0,\quad \tau^{\delta}_{k+1}
    =\inf\set{t>\tau^{\delta}_k}{\abs*{B_t-B_{\tau^{\delta}_k}}\geq \delta} 
  \end{displaymath}
  Then
  \begin{displaymath}
    \sum_{k:\tau_k^{\delta}<t}
    \abs*{A_{\tau_{k+1}^{\delta}}-A_{\tau_{k}^{\delta}}}^{p_0}
    \to c t\quad\text{in probability as $\delta\to0$}, %
  \end{displaymath}
  where $c=\E{\abs{A_{\tau^{1}_1}}^{p_0}}$.
\end{theorem}

\section{Proof of Theorem \ref{thm:2}}

Using the scaling property of the Brownian motion and of the fractional Brownian motion we have that
\begin{equation}\label{eq:scale}
  \sum_{k:\tau_k^{\delta}<t}
  \abs*{A_{\tau_{k+1}^{\delta}}-A_{\tau_{k}^{\delta}}}^{p_0}
  \eqinlaw 
  \delta^2 
  \sum_{k:\tau_k^{1}<t/\delta^2}
  \abs*{A_{\tau_{k+1}^{1}}-A_{\tau_{k}^{1}}}^{p_0}.
\end{equation}
For details see Lemma \ref{lem:scaling} below. 
Since $\E{\tau^{1}_1}=1$ we have by the law of large numbers that
\begin{displaymath}
  \lim \frac{\tau_k^1}{k}=1\quad\text{almost surely.}
\end{displaymath}
This implies that we can replace the the right hand side of \eqref{eq:scale} with
\begin{displaymath}
  \delta^2 
  \sum_{k<t/\delta^2}
  \abs*{A_{\tau_{k+1}^{1}}-A_{\tau_{k}^{1}}}^{p_0}
\end{displaymath}
and investigate its limiting behavior as $\delta\to0$.
The difference is that now the number of summands is deterministic, therefore the summands are identically 
distributed by their definition, although not independent.
It is also clear that we can further simplify the expression; to prove Theorem \ref{thm:2} it is enough to show that
\begin{equation}\label{eq:wlln}
  \frac1n
  \sum_{k<n}
  \abs*{A_{\tau_{k+1}^{1}}-A_{\tau_{k}^{1}}}^{p_0}\to\E{\abs*{A_{\tau_1^1}}^{p_0}}
  \quad\text{in probability}.
\end{equation}
As $\delta$ is now fixed to 1, we drop it from the notation. 
In this form it is a weak law of large numbers and we prove it by showing that 
the variance of this sum is $o(n^2)$. For this we use the strong mixing property 
of the increments of the fractional Brownian motion which follows easily from the 
decay of the correlation (see \cite{Maruyama1970}). We formulate it in Lemma \ref{lem:fBMmixing} below. 

We finish the proof of Theorem \ref{thm:2} by showing the next proposition.
\begin{proposition}\label{prop:1} With $B^{(k)}_t=B_{\tau_k+t}-B_{\tau_k}$, $t\geq0$, 
  and $B^{(H,k)}(x)=B^{H}(x+B_{\tau_k})-B^{H}(B_{\tau_k})$, $x\in\real$ let 
  $\xi_k=(B^{(k)}, B^{(H,k)})$. Then $(\xi_k)_{k\geq0}$ is strictly stationary and strong mixing in the sense that 
  \begin{displaymath}
    \cov(g(\xi_0),g(\xi_k))\to0 \quad \text{as} \ \ k\to\infty
  \end{displaymath}
  for any measurable $g:C[0,\infty)\times C(\real)\to \real$ functional for which $g(\xi_0)$ is 
  a square integrable random variable. 
\end{proposition}

From Proposition \ref{prop:1} we have the weak law of large numbers for squares integrables functionals $g(\xi_k)$. 
Indeed, using the stationarity 
we have the following estimation for the variance
\begin{align*}
  \Var\zjel*{\frac1n\sum_{k=0}^{n-1} g(\xi_k)}
  &=\frac1n\cov(g(\xi_0),g(\xi_0))+\frac1n \sum_{k=1}^{n-1} 2\frac{n-k}{n} \abs{\cov(g(\xi_0),g(\xi_i)}\\
  &\leq \frac2n\sum_{i=0}^{n-1} \abs{\cov(g(\xi_0),g(\xi_i)}
\end{align*}
Here $\abs{\cov(g(\xi_0),g(\xi_i)}\to0$, hence its arithmetic mean sequence does the same. 

It is possible to show that $\abs*{A_{\tau_{1}}}^{p_0}$ is square integrable, 
but we do not need this result. Indeed, if we know the $L^2$ and hence the $L^1$ 
convergence of the averages $\frac1n\sum_{k=0}^{n-1} g(\xi_k)$ for bounded $g$, 
then we have the same limiting relation 
for integrable functionals as well. 
So to finish the proof it is enough to show that $A_{\tau_{1}}$ is square integrable, 
\begin{displaymath}
  A_{\tau_1}=F(B_{\tau_1})-\int_0^{\tau_1} F'(B_s)dB_s.
\end{displaymath}
Here $F(B_{\tau_1})$ has the same law as $\int_0^1 B^H_xdx$ which has a normal law, 
so this part is obviously square intagrable. For the It\^o integral part 
we can use the izometry combined with the occupation time formula to obtain that
\begin{align*}
  \E{\zjel*{\int_0^{\tau_1} F'(B_s)dB_s}^2}
  &=\E{\int_0^\infty (F'(B_s))^2\I{s\leq\tau_1}ds}\\
  &=\int_\real \E{\zjel*{F'(x)}^2}\E{L^x_{\tau_1}}dx \\
  &\leq \sup_{\abs{x}\leq1}\E{\zjel*{B^H_x}^2}\E{\int_\real L^x_{\tau_1} dx}\\
  & = \sup_{\abs{x}\leq1}\E{\zjel*{B^H_x}^2}\E{B^2_{\tau_1}}=1.
\end{align*}

It remains to check that Proposition \ref{prop:1} holds. We do this using the 
next lemma whose proof 
involves only elementary computation, hence it is left for the reader.
\begin{lemma}\label{lem:fBMmixing}
  There is a constant depending only on $H$, such that
  for $x,x',y\in\real$ and 
  for 
  a fractional Brownian motion $B^H$ with Hurst index $H$
  \begin{displaymath}
    \cov\zjel*{\zjel{T_yB^H}_x, B^H_{x'}}\leq C\zfrac*{\abs{x}\abs{x'}}{\abs{y}^{2(1-H)}}
    ,\quad 
    \text{if} \ \ \frac{\max(\abs{x}, \abs{x'})}{\abs{y}}\leq \frac12, %
  \end{displaymath}
  \footnote{Using the Taylor expansion of $(1+x)^{2H}$ 
  \begin{align*}
    2\cov\zjel*{\zjel{T_yB^H}_x, B^H_{x'}}
    &=
    \abs{y+x}^{2H}+\abs{x'}^{2H}-\abs{y+x-x'}^{2H}
    -\zjel*{\abs{y}^{2H}+\abs{x'}^{2H}-\abs{y-x'}^{2H}}\\
    &=
    \abs{y}^{2H}\sum_{k=0}^\infty \binom{2H}{k} \zjel*{\zfrac*{x}{y}^{k}-\zfrac*{x-x'}{y}^k+\zfrac*{-x'}{y}^{k}}\\
    &=\abs{y}^{2H}\frac{xx'}{y^2}
    \sum_{k=2}^\infty \binom{2H}{k} \sum_{\ell=1}^{k-1} \binom{k}{\ell} \zfrac*{x}{y}^{\ell-1}\zfrac*{-x'}{y}^{k-\ell-1}  
  \end{align*}
  Here
  \begin{align*}
    \abs*{\sum_{k=2}^\infty \binom{2H}{k} \sum_{\ell=1}^k \binom{k}{\ell} \zfrac*{x}{y}^{\ell-1}\zfrac*{-x'}{y}^{k-\ell-1}}
    &\leq
    \sum_{k=2}^\infty \abs*{\binom{2H}{k}} \sum_{\ell=1}^{k-1} \binom{k}{\ell} \abs*{\frac{x}{y}}^{\ell-1}\abs*{\frac{x'}{y}}^{k-\ell-1}
    \\&\leq \sum_{k=2}^\infty \abs*{\binom{2H}{k}} \sum_{\ell=1}^{k-1} \binom{k-2}{\ell-1} \zfrac*12^{k-2} %
    \\&=\sum_{k=2}^\infty \abs*{\binom{2H}{k}}=\abs*{\sum_{k=2}^\infty \binom{2H}{k}(-1)^k}=\abs{1-2H}
  \end{align*}
  }
  where $T_y$ is the translation with $y$, that is $\zjel{T_y B^H}_x=B^H_{y+x}-B_y$.
\end{lemma}

This lemma extends easily with the monotone class argument to a much wider set of functionals involving scaling as well. In what follows we need it in the following form
\begin{corollary}\label{cor:1} Let $B^H$ be a fractional Brownian motion with Hurst index $H$.
  For $c>0$ let $(S_c B^H)_x=\abs{c}^{-H} B^H_{cx}$.

  Then for a measurable functional $g: C(\real)\to\real$, 
  \begin{displaymath}
    \cov\zjel*{g(S_{c_n} T_{y_n} B^{H}), g(S_{c'_n}T_{y_n'} B^H)}\to0  
  \end{displaymath}
  provided that $g(B^H)$ is square integrable and 
  \begin{displaymath}
      c_nc'_n = o((y_n-y'_n)^2).
  \end{displaymath}
\end{corollary}
\begin{proof} Using the monotone class argument it is enough to prove for functionals of the form 
  $g(B^H) = g(B^H_{x_1},\dots, B^H_{x_k})=g(B^H_{\underline{x}})$, where $g$ is a bounded continuous function 
  on $\real^k$. For this case it is enough to show that the $(S_{c_n} T_{y_n}B^H)_{\underline{x}}$ and 
  $(S_{c'_n} T_{y'_n}B^H)_{\underline{x}}$ are asymptotically independent, so eventually it is enough to 
  check that the covariances 
  $\cov((S_{c_n} T_{y_n}B^H)_{x_i},(S_{c'_n} T_{y'_n}B^H)_{x_j})$ $1\leq i,j\leq k$ 
  are vanishing in the limit.
  
  Note that 
  \begin{displaymath} 
    \cov\zjel*{(S_{c_n} T_{y_n} B^H)_{x_i},(S_{c'_n} T_{y'_n} B^H)_{x_j}}
    =(c_nc'_n)^{-H} \cov\zjel*{T_{y_n-y_n'} B^H_{c_nx_i}, B^H_{c'_nx_j}}
  \end{displaymath}
  and $\abs{c_n x_i}, \abs{c'_nx_j}<\frac12\abs{y_n-y'_n}$ for $n$ large enough.  
  Then, from Lemma \ref{lem:fBMmixing} %
  \begin{align*}
    \cov((S_{c_n} T_{y_n} B^H)_{x_i},(S_{c'_n} T_{y'_n} B^H)_{x_j}
    &=(c_nc'_n)^{-H} \cov(T_{y_n-y_n'} B^H_{c_nx_i}, B^H_{c'_nx_j})\\
    &\leq C (c_nc'_n)^{-H}\cdot \zfrac*{c_n\abs{x_i}c_n'\abs{x_j}}{\abs{y_n-y_n'}^{2(1-H)}}\\
    &\leq C\cdot \abs{x_i}\abs{x_j}\cdot\zfrac*{c_nc_n'}{(y_n-y_n')^2}^{1-H}
  \end{align*}
  which tends to zero by the assumption. The proof is complete.
\end{proof}

\begin{proof}[Proof of Proposition \ref{prop:1}]
  For the strict stationarity we need to show that $(\xi_k)_{k\geq0}$ and $(\xi_{k+\ell})_{k\geq0}$ have the same law for each $\ell\geq0$. 
  By the special structure of the sequence $(\xi_k)_{k\geq 0}$ is obtained from $\xi_0$ 
  in the same way as $(\xi_{k+\ell})_{k\geq0}$ is obtained from $\xi_\ell$.
  So it is to show that $\xi_0$ and $\xi_\ell$ has the same law. 
  It follows easily from the strong Markov property of the Brownian motion $B$ that $(\tau_\ell, B_{\tau_\ell})$ 
  is independent $B^{(\ell)}$ and $B^H$. By the stationarity of the increments of fractional Brownian motion $B^H$, 
  $(B^H_{x+y}-B^H_{y})_{x\in\real}$ has the same law as $B^H$ for any $y$. Now using the value of $B_{\tau_\ell}$ as $y$ yields 
  that $B^{(\ell)}$ and $B^{(H, \ell)}$ are independent, and has the same joint law as $B$ and $B^H$.

  For the strong mixing property it is enough to consider functionals 
  of the form $g(B,B^H)=g_1(B)g_2(B^H)$ with $g_1,g_2$ bounded and then use monotone class argument. For this special case it is enough to show that 
  \begin{displaymath}
    \E{g_1(B)g_1(B^{(k)})}\to \E^2\PEzjel*{g_1(B)}
    \quad\text{and}\quad
    \E{g_2(B^H)g_2(B^{(H,k)})\given \sigma(B)}\pto \E^2\PEzjel*{g_2(B^H)}, %
  \end{displaymath}
  as $k$ tends to infinity.
  Then using the boundedness of $g_1,g_2$ the result easily follows. 
  
  Using the translation notation from Lemma \ref{lem:fBMmixing}, 
  $\E{g_2(B^H)g_2(B^{(H,k)})\given \sigma(B)}$ is easy to express,
  \begin{displaymath}
    \E{g_2(B^H)g_2(B^{(H,k)})\given \sigma(B)}
    =\left.\E{g_2(B^H)g_2(T_xB^{H})}\right|_{x=B_{\tau_k}}.
  \end{displaymath}
  So this part follows from Corollary \ref{cor:1} and the fact that $\abs{B_{\tau_k}}\pto \infty$.

  Concerning $g_1(B)$, we can consider bounded functionals of the form $g_1(B_{\underline{t}})$, where 
  $\underline{t}\in [0,\infty)^d$ and $B_{\underline{t}}$ is the vector variable obtained from $B$ by sampling the values at the time points $\underline{t}=(t_1,\dots, t_d)$. 
  Let $\underline{t}\wedge \tau_k=(t_1\wedge \tau_k,\dots,t_d\wedge\tau_k)$, then $g_1(B_{\underline{t}\wedge \tau_k} )\pto g_1(B_{\underline{t}})$ so 
  \begin{displaymath}
    \E{g_1(B_{\underline{t}})g_1(B^{(k)}_{\underline{t}})}-\E{g_1(B_{\underline{t}\wedge\tau_k})g_1(B^{(k)}_{\underline{t}})}\to0
  \end{displaymath}
  and as $g_1(B_{\underline{t}\wedge \tau_k} )$ is independent of $g_1(B^{(k)}_{\underline{t}})$ we also have that
  \begin{displaymath}
    \E{g_1(B_{\underline{t}\wedge\tau_k})g_1(B^{(k)}_{\underline{t}})}
    =\E{g_1(B_{\underline{t}\wedge\tau_k})}\E{g_1(B^{(k)}_{\underline{t}})}
    \to\E^2\PEzjel{g_1(B_{\underline{t}})}
    \qedhere
  \end{displaymath}
\end{proof}  

\nocite{Maruyama1970}

\begin{lemma}\label{lem:scaling} For $\delta>0$ let
  \begin{equation*}
    B^{(\delta)}_t=\frac1{\delta} B_{t\delta^2},\quad
    B^{(H,\delta)}_x = \delta^{-H} B^{H}_{x\delta}
  \end{equation*}
  Then 
  \begin{enumerate}
    \item $(B,B^H)\eqinlaw (B^{(\delta)}, B^{(H,\delta)})$, 
    \item  $\tau^\delta_k(B)=\delta^2\tau^1_k(B^{(\delta)})$ and 
    \item $A_{t}(B, B^H) = \delta^{1+H} A_{t/\delta^2} (B^{(\delta)}, B^{(H,\delta)})$.
  \end{enumerate}
\end{lemma}
\begin{proof}
  The first two point follows from the scaling invariance of the (fractional) Brownian motion and from the definition of the stopping time sequence.

  For the last point
  \begin{displaymath}
    A_{t} = F(B_{t})-\int_0^{t} F'(B_u) dB_u.
  \end{displaymath}
  Here $F$ is the anti derivative of $B^H$ such that $F(0)=0$. 
  Then $B_{t}=\delta B^{(\delta)}_{t/\delta^2}$.
  If $F^{(\delta)}$ denotes the antiderivative of $B^{(H,\delta)}$ with $F^{(\delta)}(0)=0$, then 
  for positive $x$
  \begin{displaymath}
    F(x)=\int_0^x B^H_y dy = \int_0^x \delta^{H} B^{(H,\delta)}_{y/\delta}dy = \delta^{H+1} F^{(\delta)}(x/\delta),
  \end{displaymath}
  and similarly for negative $x$. From these computations $F(B_t)=\delta^{1+H} F^{(\delta)}(B^{(\delta)}_{t/\delta^2})$

  For the stochastic integral note that $F'(x)=B^H_x = \delta^{H} B^{(H,\delta)}_{x/\delta} =
  \delta^H (F^{(\delta)})'(x/\delta)$, so
  \begin{align*}
    \int_0^t F'(B_u) dB_u 
    & = \delta^H \int_0^t  (F^{(\delta)})'\zjel{\frac1{\delta} B_u} dB_u\\ %
    & = \delta^H \int_0^{t}  (F^{(\delta)})'(B^{\delta}_{u/\delta^2}) \delta dB^{(\delta)}_{u/\delta^2}\\
    & =\delta^{1+H} \int_0^{t/\delta^2} (F^{(\delta)})'(B^{\delta}_{u}) dB^{(\delta)}_{u}.
  \end{align*}
\end{proof}

\section{Proof of Theorem \ref{thm:1}}

The proof of Theorem \ref{thm:1} goes along similar lines as that of Theorem \ref{thm:2}. 
For a given interval $I=[a,b]\subset[0,\infty)$ let's define
\begin{equation}\label{eq:BIdef}
  B^{(I)}_t=\frac{B_{(b-a)t+a}-B_{a}}{(b-a)^{1/2}},\quad
  B^{(H,I)}_x = \frac{B^{H}_{(b-a)^{1/2}x+B_{a}}-B^{H}_{B_{a}}}{(b-a)^{H/2}}  
\end{equation}
The key point here is again the scaling property of the (fractional) Brownian motion:
for $0\leq s<t$, $I=[s,t]$
\begin{displaymath}
  A_{t}-A_s = F(B_t)-F(B_s)-\int_s^t F'(B_u) dB_u.
\end{displaymath}
If $F^{(I)}$ is the random function with $B^{(H, I)}$ as its derivative and $F^{(I)}(0)=0$, then 
for $x>0$
\begin{align*}
  F^{(I)}(x) &= \int_0^x B^{(H,I)}_y dy
  =\int_0^x \frac{B^{H}_{(t-s)^{1/2}y+B_{s}}-B^{H}_{B_{s}}}{(t-s)^{H/2}} dy 
  \\
  &= (t-s)^{-(H+1)/2} (F(B_s+x\sqrt{t-s})-F(B_s)) - x(t-s)^{-H/2} F'(B_s)
\end{align*}
and similarly for $x<0$. 
So if we write $B^{(I)}_{1}$ in place of $x$ we get 
\begin{displaymath}
  F^{(I)}(B^{(I)}_{1}) = (t-s)^{-(H+1)/2} \zjel*{F(B_t)-F(B_s)-(B_t-B_s)F'(B_s)}.
\end{displaymath}
Similary for the stochastic integral
\begin{align*}
  \int_0^1 (F^{(I)})'(B^{(I)}_u) dB^{(I)}_u 
  &=\int_s^t \frac{F'(B_u)-F'(B_s)}{(t-s)^{H/2}} \frac{dB_u}{(t-s)^{1/2}}\\
  &=(t-s)^{-(H+1)/2}\zjel*{\int_s^t F'(B_u) dB_u -F'(B_s)(B_t-B_s)} 
\end{align*}
So we can conclude that 
\begin{displaymath}
  \abs{A_t-A_s}^{2/(H+1)} = (t-s) g(B^{([s,t])}, B^{(H,[s,t])})
\end{displaymath}
with a suitable functional $g$.

The proof of Theorem \ref{thm:2} is based on the next two claims which we prove below.
\begin{proposition}\label{prop:2}
  For a non-degenerated interval $I\subset[0,\infty)$ %
  denote $\xi_I=(B^{(I)},B^{(H,I)})$. Then the law of $\xi_I$ does not depends on $I$ and
  \begin{displaymath}
    \cov(g(\xi_{I_n}), g(\xi_{J_n}))\to 0 \quad 
    \text{if} \ \ \max(\abs{I_n},\abs{J_n})\to0 \ \ \text{and} \ \ \liminf_n \dist(I_n,J_n)>0,
  \end{displaymath} 
  provided that $g(\xi_{[0,1]})$ is square integrable.
\end{proposition}

\begin{corollary}\label{cor:wlln}
  Suppose that $(\pi_n)$ is a sequence of subdivisions of $[0,t]$, such that the mesh
  $\max_{I\in\pi_n}\abs{I}\to0$ and $g$ is a functional such that
  $g(\xi_{[0,1]})$ is integrable. Then 
  \begin{displaymath}
    \sum_{I\in\pi_n} \abs{I}g(\xi_I)\to t\E{g(\xi_{[0,1]})}\quad\text{in $L^1$}
  \end{displaymath}
\end{corollary}

Theorem \ref{thm:2} follows from  %
Corollary \ref{cor:wlln} if we apply it to 
\begin{displaymath}
  g(\xi_{[0,1]})= \abs*{F(B_1)-\int_0^1 F'(B_s)dB_s}^{p_0},\quad F'(x)=B^H_x, \quad F(0)=0.
\end{displaymath}
We check that $g(\xi)$ is integrable by showing that $g(\xi_{[0,1]})\in L^{2/p_0}$. Note that
$2/p_0<2$ so it is enough to show that the next random variable is square integrable
\begin{displaymath}
  F(B_1)-\int_0^1 F'(B_s)dB_s,\quad F'(x)=B^H_x, \quad F(0)=0.
\end{displaymath}
Since by trivial estimations $\E{F(x)^2}<x^4$ we obviously 
have that $F(B_1)\in L^2$. For the second term we apply again \ito/ izometry 
followed by the occupational time formula to get that
\begin{align*}
  \E{\zjel*{\int_0^1 F'(B_s)dB_s}^2}
  &=\int_{\real} \E{\zjel{B^H_x}^2} \E{L^x_1} dx\\
  &=\int_{\real} \abs{x}^{2H} \E{L^x_1}dx\\
  &\leq \int_\real (x^2+1)\E{L^x_1}dx 
  =\E{\int_0^1 \zjel*{B_s^2+1} ds}=3/2. \qedhere
\end{align*}

\begin{proof}[Proof of Proposition \ref{prop:2}]
  We start with the law of  $\xi_I$.
  Since the increments of the fractional Brownian motion $B^H$ are stationary, 
  $(B^H_{x+y}-B^H_y)_{x\in \real}$ has the same law as $B^H$ for any deterministic $y$.
  Then by the independence of $B$ and $B^{H}$ the conditional law of 
  $(B^{H}_{x+B_a}-B^H_{B_a})_{x\in\real}$ given $B$ does 
  not depend on $B$, that is $(B^{H}_{x+B_a}-B_{B_a})_{x\in\real}$ is independent of $B$ with 
  the same law as $B^H$. But then by the scaling invariance of $B^H$ the same is true for $B^{(H,I)}$. 
  Finally, by the Markov property and scaling invariance of $B$ we get that $B^{(I)}$ is also
  a Brownian motion which is obviously independent of $B^{(H,I)}$. So $\xi_I$ has the same law as $(B,B^H)$ 
  which is $\xi_{[0,1]}$ by definition.

  To show the asymptotics of the covariance it is enough to consider again functionals of the form 
  $g(B,B^H)=g_1(B)g_2(B^H)$ %
  where $g_1$, $g_2$ are bounded.  As in the proof of Proposition \ref{prop:1} it is enough to show that
  $\E{g_2(B^{(H,I_n)})g_2(B^{(H,J_n)})\given \sigma(B)}\pto \E^2\PEzjel*{g_2(B^{H})}$ and that 
  $\E{g_1(B^{(I_n)}g_1(B^{(J_n)})}\to \E^2\PEzjel*{g_1(B)}$ whenever
   $(I_n=[a_n,b_n],J_n=[c_n,d_n])_{n\geq1}$ is a sequence such that 
   $\max(\abs{I_n},\abs{J_n})\to0$ and $\inf_n\dist(I_n,J_n)>0$.

  Using the independence of $B$ and $B^H$ we get that
  \begin{displaymath}
    \E{g_2(B^{(H,I_n)})g_2(B^{(H,J_n)})\given \sigma(B)}=
    \left.\E{g_2(S_{\abs{I_n}^{1/2}} T_x B^{H}) g_2(S_{\abs{J_n}^{1/2}} T_y B^{H})}\right|_{x=B_{a_n}, y=B_{c_n}},
  \end{displaymath}
  where  we used the notation of Corollary \ref{cor:1}. By assumption 
  \begin{displaymath}
    \frac{\abs{I_n}^{1/2}\abs{J_n}^{1/2}}{(B_{c_n}-B_{a_n})^2}\pto 0 \quad\text{as $n\to\infty$}.
  \end{displaymath}
  Then Corollary \ref{cor:1} shows that $\E{g_2(B^{(H,I_n)})g_2(B^{(H,J_n)})\given \sigma(B)}\pto\E^2\PEzjel*{g_2(B^H)}$.

  For $\E{g_1(B^{(I_n)}) g_1(B^{(J_n)})}=\E{g_1(S_{\abs{I_n}} T_{a_n}B))g_1(S_{\abs{J_n}} T_{c_n}B)}$ we can also use Corollary \ref{cor:1} with $H=1/2$.
\end{proof}

\begin{proof}[Proof of Corollary \ref{cor:wlln}]
  First assume that $g(\xi_{[0,1]})$ is square integrable, 
  and for a sequence of partitions $\pi_n$ of $[0,t]$ let 
  \begin{displaymath}
    f_n = \sum_{I,J\in \pi_n} \I[I\times J] \cov(g(\xi_{I}), g(\xi_J))
  \end{displaymath}
  Then 
  \begin{displaymath}
    \Var \zjel*{\sum\nolimits_{I\in\pi_n} \abs{I}g(\xi_I)}=\int_{[0,t]^2} f_n \to 0,
    \quad\text{if $\max_{I\in\pi_n}\abs{I_n}\to0$,}
  \end{displaymath}
  since the sequence of function $(f_n)$ is dominated by $\Var g(\xi_{[0,1]})$ and tends 
  to zero everywhere but the diagonal of $[0,t]^2$ by Proposition \ref{prop:2}.
  Also by Proposition \ref{prop:2} the expectation $\E{\sum_{I\in\pi_n} \abs{I}g(\xi_I)}=t\E{g(\xi_{[0,1]})}$
  does not depend on $n$, hence we have the claim for $g$ provided that $g(\xi_{[0,1]})\in L^2$.

  For general $g$, when $g(\xi_{[0,1]})$ is integrable, we can use truncation. 
\end{proof}

\section{The zero energy part of the median process}

Let $B$ be a Brownian motion and suppose that $\zjel*{D_t(x)}_{t \geq 0, x \in [0,1]}$
satisfies the following stochastic differential equation
\begin{equation}\label{eq:D}
  d D_t (x)=\sigma (D_t(x))dB_t, 
  \quad \sigma(x)=x\wedge(1-x), \quad D_0(x)=x, \quad x\in[0,1].
\end{equation}
This two parameter process was analyzed in \cite{Prokajetal2011} in detail.
It was shown that the solution $D_t(x)$ can be viewed as a conditional distribution 
function, and this justifies the (conditional) quantile name for the process 
$q_t=D^{-1}\zjel*{\alpha}$, $\alpha\in(0,1)$, and particularly 
the (conditional) median name for the process $m_t=D^{-1}\zjel*{\frac{1}{2}}$.

In \cite{ProkajBondici2022} it was proved that the quantile process 
$q_t$ is not a semimartingale, so neither the median process $m_t$ is it.
The following decomposition formula holds for $q_t$ 
(cf. Subsection 5.2 in \cite{ProkajBondici2022})
\begin{displaymath}
  A_t = q_t + \int_0^t \zjel*{D_s^{-1}}'(\alpha)\sigma(\alpha)dB_s.
\end{displaymath}
This $A$ is a process of zero energy, that is, the quadratic variation of $A$ 
exists and $\br{A}\equiv 0$. If $A$ would have finite total variation, then 
$q_t$ would be a semimartingale, so $A$ should have infinite total variation.

We will refer to the process $A$ as the zero energy part
of the quantile process. In the following we prove that the local
martingale part in the previous decomposition of $q_t$ is a true martingale
and we will investigate the following main question: however
the total variation of $A$ is infinite, and the quadratic variation is identically $0$, 
is there any $p\in (1,2)$ for which the $p$-variation of $A$ is positive and finite?
We are not able to give a mathematically rigorous answer to this latter question, but
we formulate some heuristic arguments which are supported by some simulation results.

\subsection{Space inverse of a stochastic flow}

In this subsection we revise a method for obtaining the space inverse
of a stochastic flow (which is given by an \ito/ diffusion) at a given time point.
We will use this in the next subsection (in the proof of the martingale 
property of the local martingale part of $q_t$).
We prove only for the case of the unit diffusion coefficient but with 
suitable transformations this result can be extended.

Let $B$ be a Brownian motion and 
let $(G_t(x))_{t\geq 0, x\in \real}$ 
be a stochastic flow which satisfies the following equation
\begin{equation}\label{eq:G_gen}
  dG_t(x) = \mu(G_t(x))dt + dB_t, \quad G_0(x) = x, \quad x\in\real,
\end{equation}
where $\mu$ is a bounded measurable function.
Suppose that on an almost sure event the mapping $(t,x)\mapsto G_t(x)$
is continuous and the mappings 
$x\mapsto G_t(x)$ are homeomorphisms of $\real$ for all $t$. 
We now define a process which produces as its terminal 
value on an almost sure event the space inverse of $G_T(x)$ for a certain $T>0$. 
For this we will use the time reversion of the Brownian motion.

For a fixed time horizon $T>0$ let $\tB^T_u = B_{T-u}-B_T$ 
be the time reversed Brownian motion on $[0,T]$. 
Let $(H_t^T(x))_{t\in [0,T], x\in \real}$ 
be the solution of the following equation
\begin{equation}\label{eq:H_gen}
  dH_t^T(x) = -\mu(H_t^T(x))dt + d\tB^T_t, \quad H_0^T(x)=x, \quad x\in\real.
\end{equation}
Suppose again that on an almost sure event the mapping 
$(t,x)\mapsto H^T_t(x)$ is continuous and the mappings 
$x\mapsto H^T_t(x)$ are homeomorphisms of $\real$ for all $t$. 

Recall the following result from \cite{Davie2007} (see also in \cite{Flandoli2011}; 
a different approach and a generalization of this result can be found in
\cite{Shaposhnikov2016} and \cite{Shaposhnikov2017})
\begin{theorem}
  Consider the following equation
  \begin{equation}\label{eq:int}
    X_t = X_0+B_t+\int_0^t f(s,X_s)ds, \quad t\geq 0,
  \end{equation}
  where $X_t\in\real^d$, $B$ is a standard 
  $d$-dimensional Brownian motion and $f$ is a bounded 
  Borel function from $[0,\infty)\times \real^d$ to $\real^d$.
  
  For almost every Brownian path $B$, there is a unique continuous
  $X:[0,\infty)\to \real^d$ satisfying \eqref{eq:int}.
\end{theorem}
It will be pointed out in the proof of the next proposition the role of this theorem.
Applying this result, we know that on an almost sure event there 
are unique continuous  
$\hat{G}(x_0):[0,\infty)\to\real$ and $\hat{H}(x_0):[0,\infty)\to\real$
which satisfy
\begin{displaymath}
  d\hat{G}_t(x_0) = \mu(\hat{G}_t(x_0))dt + dB_t, \quad \hat{G}_0(x_0) = x_0,
\end{displaymath}
and
\begin{displaymath}
  d\hat{H}_t(x_0) = -\mu(\hat{H}_t(x_0))dt + dB_t, \quad \hat{H}_0(x_0) = x_0.
\end{displaymath}
Taking into account all rational $x_0$ and then using the continuity of 
$G_t(x)$ and $H^T_t(x)$, 
we can conclude that on an almost sure event $\Omega' \subset \Omega$ 
the mappings $(t,x)\mapsto G_t(x)$ and $(t,x)\mapsto H_t(x)$ are unique.
Then we have the following
\begin{proposition}\label{prop:space inv of flow}
  On an almost sure event, for all $u \in [0,T]$ and for all $x \in\real$ 
  we have the following: 
  $G_{T-u}(x) = H^T_u(G_T(x))$ and $H^T_{T-u}(x) = G_u(H^T_T(x))$. 
  
  Especially, we have $H^T_T(G_T(x))=x$ and $G_T(H^T_T(x))=x$, 
  so $G^{-1}_T = H_T^T$.  
\end{proposition}
\begin{proof}
  We restrict ourselves to $\Omega'$. 
  Let $0 \leq s < t \leq T$. For $u \in [0,T]$ let $\tilde{G}_u^T(x) =
  G_{T-u}(x)$ be the time reversed process of $\zjel*{G(x)}_{t\in[0,T]}$. 
  Using \eqref{eq:G_gen} we obtain
  \begin{multline*}
    \tilde{G}_t^T(x) - \tilde{G}_s^T(x) = -(G_{T-s}(x)-G_{T-t}(x))
    = -\int_{T-t}^{T-s}\mu(G_u(x))du -(B_{T-s}-B_{T-t})=\\ 
    -\int_s^t \mu\zjel*{G_{T-u}(x)}du + \zjel*{B_{T-t}-B_T-(B_{T-s}-B_T)}
    = -\int_s^t \mu\zjel*{\tilde{G}_u^T(x)}du + \zjel*{\tB_t^T - \tB_s^T}.
  \end{multline*}
  Substituting $G^{-1}_T(x)$ in place of $x$ yields
  \begin{displaymath}
    \tilde{G}_t^T(G^{-1}_T(x)) - \tilde{G}_s^T(G^{-1}_T(x)) = 
    -\int_s^t\mu\zjel*{\tilde{G}_u^T(G^{-1}_T(x))}du 
    + \zjel*{\tB_t^T - \tB_s^T}.
  \end{displaymath}
  As $\tilde{G}_0^T\zjel*{G^{-1}_T(x)} = G_T\zjel*{G^{-1}_T(x)} = x$, 
  we have obtained that
  $\tilde{G}_u^T\zjel*{G^{-1}_T(x)}$ solves the equation \eqref{eq:H_gen} of
  $H^T(x)$. 
  However, from \cite{Veretennikov1980} we know that \eqref{eq:H_gen} 
  has a unique strong solution, but it is not obvious that 
  $\zjel*{\tilde{G}_u^T\zjel*{G^{-1}_T(x)}}_{u\in [0,T]}$ is adapted to the filtration 
  generated by $\tB^T$, hence we need the previously cited uniqueness result. 
  It follows that on $\Omega'$ we have
  $H^T_u(x) = \tilde{G}_u^T\zjel*{G^{-1}_T(x)}$, so $H^T_u\zjel*{G_T(x)} = G_{T-u}(x)$.
  
  A similar argument for $H^T$ yields $G_u\zjel*{H^T_T(x)}=H^T_{T-u}(x)$.
\end{proof}

\subsection{
The local martingale part of the quantile process is a true martingale}

Now we turn to the case of the quantile process.
Using the result of the previous subsection we derive a process which produces 
a transformed version of $D_T^{-1}$ as its terminal value. 
In the following we do not use the exact form of $\sigma$ 
from the equation \eqref{eq:D} of $D(x)$, we only use that it is 
a Lipschitz continuous function with Lipschitz constant $L$, 
but we also provide the particular 
results for that special case.

As the first step, to be able to apply Proposition 
\ref{prop:space inv of flow}, consider the Lamperti transform of $D(x)$, 
so let $\psi:C\to\real$ be such that $\psi'=1/\sigma$, where
$C$ is a connected component of $\real\setminus\event{\sigma=0}$.
In the special case of $\sigma(x)=x\wedge(1-x)$ and $C=(0,1)$ 
we can use $\psi(x)= \sign(1-2x)\ln(1-\abs{1-2x})$.
Let $G_t(x) = \psi\circ D_t \circ\psi^{-1}(x)$. Then $G(x)$ evolves as
\begin{equation}\label{eq:G}
  dG_t(x)=-\frac12(\sigma'\circ\psi^{-1})\zjel*{G_t(x)}dt + dB_t, 
  \quad G_0(x)=x, \quad x\in\real.
\end{equation}
Let $H_t^T(x)$ be the solution of 
\begin{equation}\label{eq:H}
  dH_t^T(x)=\frac12(\sigma'\circ\psi^{-1})\zjel*{H_t^T(x)}dt+d\tB_t^T, 
  \quad H_0^T(x)=x, \quad x\in\real.
\end{equation}
In the special case the equations are
\begin{equation}
  dG_t(x)=\frac12\sign(G_t(x))dt + dB_t,
  \quad G_0(x)=x, \quad x\in\real, \tag{\ref*{eq:G}${}^{\prime}$}
\end{equation}
and
\begin{equation}
  dH_t^T(x)=-\frac12\sign(H_t^T(x))dt+d\tB_t^T, 
  \quad H_0^T(x)=x, \quad x\in\real. \tag{\ref*{eq:H}${}^{\prime}$}
\end{equation}

In order to be able to calculate $\zjel*{D_T^{-1}}'$ easily, 
let $s$ be a scale function (or Zvonkin transform, \cite{Zvonkin1974}) 
for $H^T(x)$. This transform removes the drift, and $s$ satisfies 
\begin{displaymath}
  s'(x)(\sigma'\circ\psi^{-1})(x) + s''(x)=0.
\end{displaymath}
For such a function we have $\displaystyle s'=\frac{c}{\sigma\circ\psi^{-1}}$ 
with some $c\in\real$, and the transformed process 
$\zjel*{F_t^T(x)=s\circ H_t^T(x)}_{t\geq 0, x\in\real}$ satisfies
\begin{equation}\label{eq:F_gen}
  dF_t^T(x)=\frac{c}{\sigma\circ(s\circ\psi)^{-1}} \zjel*{F_t^T(x)}d\tB_t^T, 
  \quad F_0^T(x)=x, \quad x\in\real.
\end{equation}
In the special case of the quantile process a possible $s$ can be 
$s(x) = \sign(x)(\exp{\abs{x}}-1)$, and with this choice
\begin{equation}\label{eq:F}
  dF^T_t(x) = \zjel*{1+\abs*{F^T_t(x)}}d\tB^T_t, \quad F^T_0(x) = x, \quad x\in\real. 
  \tag{\ref*{eq:F_gen}${}^{\prime}$}
\end{equation}
To calculate $(F_T^T)'=f_T^T$, we can use Lemma 16 from 
\cite{Prokajetal2011}. Since 
\begin{displaymath}
  \zjel*{\frac{c}{\sigma\circ(s\circ\psi)^{-1}}}'=
  -\sigma'\circ(s\circ\psi)^{-1},
\end{displaymath}
we know that $\displaystyle \frac{c}{\sigma\circ(s\circ\psi)^{-1}}$ 
is a Lipschitz continuous function (with the same 
Lipschitz constant $L$ as $\sigma$), 
so $F_t^T$ is differentiable in its space variable, 
and the space derivative $f_t^T$ satisfies  
\begin{displaymath}
  df_t^T(x) = -\sigma'\circ(s\circ\psi)^{-1}(F_t^T(x))f_t^T(x)d\tB_t^T, 
  \quad f_0^T(x)=1, \quad x\in\real.
\end{displaymath}
From this we obtain that
\begin{equation}\label{eq:f_t^T explicit}
  f_t^T(x)=\exp{N_t^T(x)-\frac12[N^T(x)]_t},
  \quad
  N_t^T(x) = -\int_0^t \sigma'\circ(s\circ\psi)^{-1}(F_u^T(x))d\tB_u^T
\end{equation}
and this yields the following estimation
\begin{displaymath}
  \E{\abs*{f_t^T}^p} \leq 
  \exp{L^2 p\zjel*{p-\frac12} t},
\end{displaymath}
where $L$ is the Lipschitz constant for $\sigma$.

To be able to apply Proposition \ref{prop:space inv of flow}, we need to guarantee
that on an almost sure event $(t,x)\mapsto G_t(x)$ and $(t,x)\mapsto H^T_t(x)$ 
are continuous, $G_t(x)$ and $H^T_t(x)$ are homeomorphisms of $\real$ for all $t$
and its drifts are bounded. The latter property follows from the Lipschitz continuity of $\sigma$.
Moreover, since also $\displaystyle \frac{c}{\sigma\circ(s\circ\psi)^{-1}}$
is a Lipschitz con\-ti\-nuous function, we know that 
(Theorems 37 and 46 from Chapter V in \cite{Protter2005}) 
the above properties hold for 
$\zjel*{D_t(x)}_{t\geq 0, x\in\real}$ and $\zjel*{F^T_t(x)}_{t\geq 0, x\in\real}$,
and hence hold for $\zjel*{G_t(x)}_{t\geq 0, x\in\real}$ and $\zjel*{H^T_t(x)}_{t\geq 0, x\in\real}$
(as they are transformations).
So by Proposition \ref{prop:space inv of flow} we know that almost surely
$H_T^T = G_T^{-1} = \psi\circ D_T^{-1}\circ\psi^{-1}$, so  
\begin{displaymath}
  D_T^{-1} = \zjel*{s\circ\psi}^{-1}\circ F_T^T\circ\psi,
\end{displaymath}
hence 
\begin{displaymath}
  \zjel*{D_T^{-1}}' = 
  \frac{c\cdot \sigma^2\circ(s\circ\psi)^{-1}\circ F_T^T \circ \psi 
  \cdot f_T^T \circ \psi}{\sigma}.   
\end{displaymath}
Using these results we can prove that the local martingale part of $q_t$ 
is a true martingale.
\begin{lemma}
  If $\abs*{\sigma}\leq K$ on $C$, then
  \begin{displaymath}
    \sup_{0\leq u\leq t}\E{\zjel*{D_u^{-1}}'(\alpha)}^2 < \infty
  \end{displaymath}
  for every $t\geq 0$ and $\alpha\in C$.
\end{lemma}
\begin{proof}
  Using the previous results we have
  \begin{align*}
    \E{\zjel*{\zjel*{D_u^{-1}}'(\alpha)\sigma(\alpha)}^2}
    &=c^2 \E{\sigma^4\circ(s\circ\psi)^{-1}\circ F_u^u \circ \psi(\alpha) 
    \cdot f_u^u \circ \psi(\alpha)} \\
    &\leq 
    c^2 
    \zjel*{\E{\sigma^8\circ(s\circ\psi)^{-1}\circ F_u^u \circ \psi(\alpha)}}^{1/2}
    \zjel*{\E{\zjel*{f_u^u}^2 \circ \psi(\alpha)}}^{1/2}\\
    &\leq
    c^2 K^4 \exp{\frac32 L^2 u}, 
  \end{align*}
  hence
  \begin{displaymath}
    \sup_{0\leq u\leq t}\E{\zjel*{D_u^{-1}}'(\alpha)}^2 \leq
    \frac{c^2 K^4}{\sigma^2(\alpha)} \exp{\frac32 L^2 t}. 
  \end{displaymath}
\end{proof}

\begin{corollary}
  In the case of $\sigma(x)=x\wedge(1-x)$ and $C=(0,1)$ it follows that 
  the local martingale part of the quantile process is a 
  true martingale.
\end{corollary}

\subsection{Simulation framework and results}

In this subsection we will restrict ourselves to the case 
$\alpha=\frac12$, so $q_t=m_t$ is the median process, and with the previous 
notations we have 
\begin{equation}\label{eq:medfromF}
  m_T = (s \circ \psi)^{-1}\zjel*{F_T^T(0)}.  
\end{equation}
In the next we present the simulation framework and the results which supports 
the following
\begin{conj} With the above notations, for $p=\frac43$, the stochastic limit of
  \begin{equation}\label{eq:pvar A}
    \sum_{k=0}^{n-1}\abs*{A_{t_{k+1}}-A_{t_{k}}}^p, \quad n\to \infty 
  \end{equation}
  is positive and finite (where $0=t_0 < t_1 < \dots < t_n=T$ denotes, 
  with simplified indices, a sequence of partitions of $[0,T]$ 
  which have grid size tending to $0$).
\end{conj}

Instead of the sum in the above conjecture, we investigate the following sum
\begin{equation}\label{eq:pvar A cond}
  \sum_{k=0}^{n-1}\abs*{\E{A_{t_{k+1}}-A_{t_{k}}\given \F_{t_k}}}^p,
\end{equation}
where $\F$ is the natural filtration of the Brownian motion $B$. 
If for $p>1$ the $p$-adic variation \eqref{eq:pvar A} tends to $0$ in expectation, 
then so does the conditional version \eqref{eq:pvar A cond}, 
so if the conditional version has a positive and 
finite limit, then \eqref{eq:pvar A} can not tend to $0$.
Moreover,
in this way the martingale term in the definition of $A$ can be eliminated
\begin{displaymath}
  \E{A_{t_{k+1}}-A_{t_{k}}\given \F_{t_k}} = \E{m_{t_{k+1}}-m_{t_{k}}\given \F_{t_k}} =
  \E{m_{t_{k+1}}\given \F_{t_k}} - m_{t_k}.
\end{displaymath}

\subsubsection{Calculating the median and calculating the conditional expectation}

For calculating the value of the median we can use \eqref{eq:medfromF}, since for $m_T$
it is enough to calculate $F^T_T(0)$ (and transform its value). 
In order to do this, we want to use the
discretised version of the equation \eqref{eq:F}
\begin{equation}\label{eq:discrF}
  F^T_{t_{i+1}}=F^T_{t_{i}}+\zjel*{1+\abs*{F^T_{t_{i}}}}\Delta \tB^T_{t_i},
  \quad F_{t_0}=0,
\end{equation}
where $0=t_0<t_1<\dots<t_k=T$ is an equidistant grid of $[0,T]$ with mesh size 
$\Delta t = \frac{T}{k}$,
and $\Delta \tB_{t_i}^T$ are correspondingly rescaled independent Rademacher variables 
(with expectation $0$ and variance $\Delta t$, so 
$\P\zjel*{\Delta \tB_{t_i}^T=\sqrt{\Delta t}}=
\P\zjel*{\Delta \tB_{t_i}^T=-\sqrt{\Delta t}}=
\frac{1}{2}$).
We restrict ourselves to the interval $[0,1]$, so $T\in [0,1]$. For the calculation of the 
median we use the sequence $(t_i=\frac{i}{n})_{i=0,\dots,k}$ for some 
suitable values of $n$ (typically powers of $10$). 

Between the time reversed Brownian motions we have the following relationship: 
for $T_1 < T_2$ and $u \in [0, T_2]$, we have
\begin{equation}\label{eq:rel BMs}
  \tB_u^{T_2} = \tB^{T_1}_{(u-(T_2 - T_1))\vee 0}+\tB^{T_2}_{u\wedge (T_2-T_1)}, 
\end{equation}
which easily follows from the definition of $\tB^T$. Since $\tB^{T_1}_{(u-(T_2 - T_1))\vee 0}$
can be calculated from the increments of $B$ in the interval $[0,T_1]$,
and $\tB^{T_2}_{u\wedge (T_2-T_1)}$ is an increment of $B$ in $[T_1,T_2]$, they are independent.

Suppose that we want to calculate the values of the median in two consecutive time points, 
in $T_1$ and $T_2=T_1+\Delta t$, on the same trajectory (until $T_1$). 
Suppose that we use a sequence $(\Delta\tB^{T_1})$ for $m_{T_1}$; 
in order to calculate $m_{T_2}$, we append one independent rescaled Rademacher variable 
to the beginning of the sequence $(\Delta\tB^{T_1})$, and we calculate $m_{T_2}$ along 
this extended sequence. 
In this way we can guarantee that we remain on the same trajectory.

Now we describe how to estimate the conditional expectation of $m$. Let $T_1$ and 
$T_2=T_1+\Delta t$ be again two consecutive time points. 
We want to calculate one realization of $\E{m_{T_1+\Delta t}\given \F_{T_1}}$
on the same trajectory as we did in the case of $m_{T_1}$. 

Since for $T_1<T_2$ and $t \in [T_2-T_1,T_2]$ we have 
$F^{T_2}_t(x)=F^{T_1}_{t-(T_2-T_1)}\zjel*{F^{T_2}_{T_2-T_1}(x)}$,
we obtain that $F^{T_2}_{T_2}(0)=F^{T_1}_{T_1}\zjel*{F^{T_2}_{\Delta t}(0)}$, 
and from \eqref{eq:rel BMs} we know that $F^{T_2}_{\Delta t}(0)$ is independent of 
$F^{T_1}_{T_1}(x)$.

Suppose that $m_{T_1}$ is calculated using the sequence $(\Delta \tB^{T_1})$. 
To calculate the conditional expectation $\E{m_{T_1+\Delta t}\given \F_{T_1}}$ along
this trajectory, we can do the following. We use the discretised equation \eqref{eq:discrF}
of $F$ with the same driving sequence $\Delta \tB^{T_1}$, but with different initial
values: we approximate the distribution of $F^{T_2}_{\Delta t}(0)$ using a finer grid of 
$[0,\Delta t]$ (using \eqref{eq:discrF} of course, with independent random increments), 
and we calculate the values of $F^{T_1}_{T_1}$ from these points as initial values. 
Finally, we take the average of $F^{T_1}_{T_1}$'s.

\subsubsection{One single increment and \texorpdfstring{$p$}{p}-variation for some values of \texorpdfstring{$p$}{p}}

Next we summarize the simulation results. First, we present the results regarding one single increment 
on the time interval $[1,1+\Delta \tau]$ for different values of $p$, 
so we investigate 
\begin{displaymath}
  \abs*{\E{A_{1+\Delta \tau}-A_1 \given \F_1}}^p = 
  \abs*{\E{m_{1+\Delta \tau}\given \F_1}-m_1}^p.
\end{displaymath}
The values of $\Delta \tau$ range from $10^{-2}$ to $10^{-7}$. 
To approximate the values of $m_1$ and $\E{m_{1+\Delta\tau}\given \F_1}$
we set the grid mesh to $10^{-6}$ ($n=10^6$ gridpoints). 

We have simulated $10^3$ realizations, and the results can be seen on the figure 
below. This is a log-log (base $10$) plot. 
On the $x$-axis we indicated $\Delta\tau$,
while on the $y$-axis we marked the average of the above mentioned increments. 
Around the mean we can also see the $95\%$ confidence interval of the expected value.
We also give the slope and intercept values for the fitted lines, and the sum of 
squared residuals (SSR).

\begin{table}[h]
  \centering
  \begin{minipage}[]{0.40\textwidth}
    \includegraphics[width=1.0\textwidth]{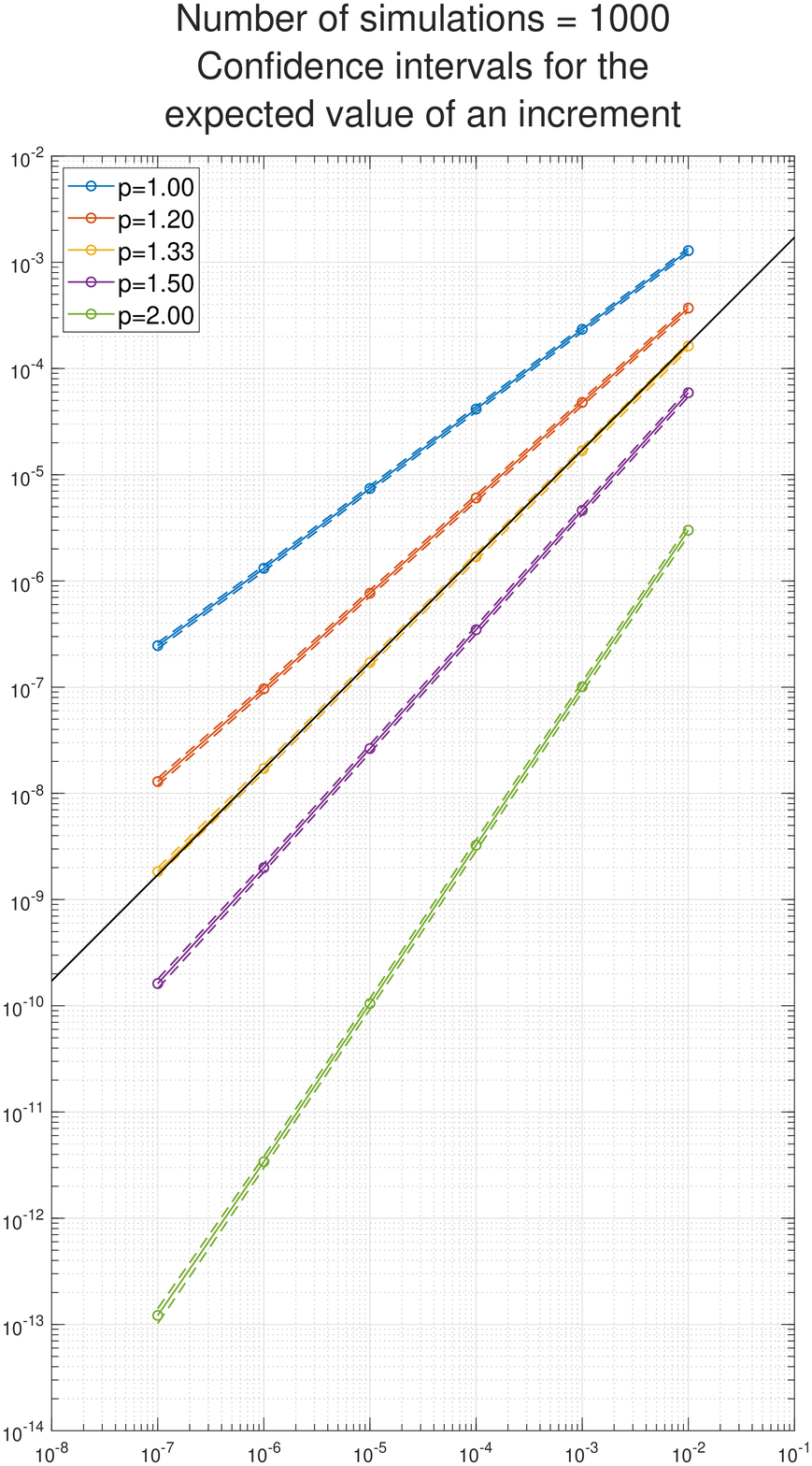}
    \captionof{figure}{One single increment}
    \label{fig:1}
  \end{minipage}
  \hfill
  \begin{minipage}[]{0.50\textwidth}
    \centering
    \begin{tabular}{|c|c|c|c|}
      \hline
      p & slope & intercept & SSR\\
      \hhline{|=|=|=|=|}
      1.00 & 0.7455 & -1.3994 & 0.0002\\ 
      1.20 & 0.8935 & -1.6443 & 0.0002\\ 
      1.33 & 0.9920 & -1.8034 & 0.0003\\ 
      1.50 & 1.1149 & -1.9980 & 0.0004\\ 
      2.00 & 1.4820 & -2.5593 & 0.0008\\ 
      \hline
    \end{tabular}
    \captionof{table}{Fitted lines}
    \label{tab:1}
  \end{minipage}
\end{table}

This figure (Figure \ref{fig:1}) and the fitted lines (Table \ref{tab:1}) 
suggest the following: there is a linear relationship between the logarithm
of $\Delta\tau$ and the logarithm of the expected value of one single increment: 
\begin{equation*}
  \lg\zjel*{\E{\abs*{\E{A_{1+\Delta \tau}-A_1 \given \F_1}}^p}}=c_0+c_1\lg\zjel*{\Delta \tau}
  \; \Rightarrow \;
  \E{\abs*{\E{A_{1+\Delta \tau}-A_1 \given \F_1}}^p}=\tc_0 \zjel*{\Delta\tau}^{c_1}.
\end{equation*}
By a scaling argument we can suppose that this relationship is valid not only on 
the interval $[1,1+\Delta\tau]$ but also on other intervals.
The slope of the thick line is $c_1 = 1$, and is very close to the points which belong to 
the case $p=\frac43$. As the number of intervals in $[0,1]$ with length $\Delta\tau$ 
is roughly $\frac{1}{\Delta\tau}$, this has the following consequences: 
\begin{itemize}
  \item for $p=\frac43$, the sum
  $\sum_{k=0}^{n-1}\E{\abs*{\E{A_{t_{k+1}}-A_{t_k} \given \F_{t_k}}}^p}$
  has a positive and finite limit;
  \item for $p<\frac43$, the above sum is unbounded from above;
  \item for $p>\frac43$, the above sum tends to $0$.
\end{itemize}

Next, we present some results regarding 
the $p$-variation processes for different values of $p$.
On the figure below (Figure \ref{fig:2}) we investigate the behavior of the $p$-variations
at time $t=1$ . 
We have simulated $500$ trajectories with different step values
(the mesh $\Delta t$ of the grid ranges from $10^{-2}$ to $10^{-5}$, which 
corresponds to $n=10^2, \dots ,10^5$ gridpoints).

The results are in line with the previous ones.
We can observe that in the case of $p=\frac43$ the distributions of the $p$-variation 
at $t=1$ get more contcentrated around the same value as $\Delta t$ decreases, while
for $p<\frac{4}{3}$ and $p>\frac43$ they get contcentrated around increasing, respectively
decreasing values.
\begin{center}
  {
  \includegraphics[width=0.60\textwidth,keepaspectratio]{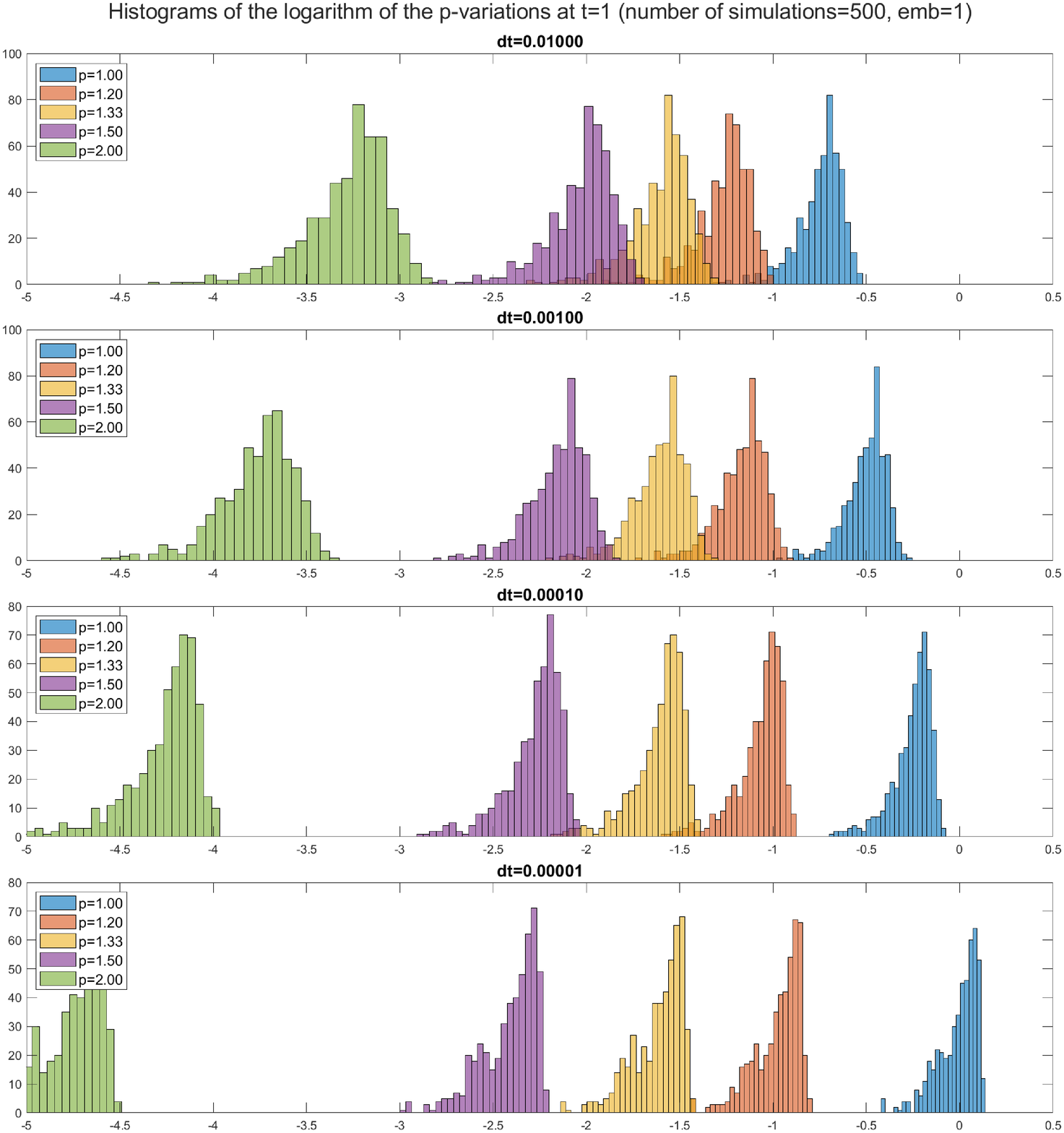}  
  \captionof{figure}{}
  \label{fig:2}
  }
\end{center}


\end{document}